# Traffic Management in Smart Cities Using the Weighted Least Squares Method

Hazim Al Gharrawi [1] and Majid Bani Yaghoub[2]

Department of Civil and Mechanical Engineering
University of Missouri-Kansas City
5110 Rockhill Rd Kansas City, MO, 64110

[2]Department of Mathematics and Statistics
University of Missouri-Kansas City
Kansas City, Missouri 64110-2499, USA

## Abstract

The Weighted Least Squares (WLS) method has several applications in science, engineering, and technology. In this paper we apply the WLS method to find the best location of a portable drone station for traffic management services. As a case study, we show that the optimal location of a drone station bounded by a polygon in Kansas City, Missouri can be determined by applying the WLS method to daily traffic congestion data. Given the mobility of the drone station, the optimal location can be updated based on the morning and afternoon rush hour data. This paper highlights the importance of WLS optimization methods in smart city planning.

## 1. Introduction

A smart city is an urban area with advanced technology for sustainable roads, optimized resource consumptions, digital government and well-connected citizens. Software reader technology such as license plate reading, facial recognition system can used as a surveillance system in smart cities to increase safety of citizens [1]. In addition, unmanned air vehicles such as drones are increasingly used in several civilian situations. A drone often communicates through wireless with a ground control station, and it may be operated remotely by human or autonomously by software. In the recent years, drones have been used for traffic surveillance and management [2].

Recent developments and technologies applied to smart cities have enabled researchers to use mathematical and computational methods to promote smart cities [3]. This includes computational urban planning and management [4], intelligent transportation systems [5], and prediction of wide traffic congestion in smart cities.

In this paper we demonstrate how one of the most powerful methods in mathematics and statistics can be used to optimally control and reduce traffics in smart cities. The main goal of the present work is to use the least squares method [5] to identify the best location of a portable drone station in a metropolitan city. By the best location we mean a location that is closest to traffic congestions and accidents that may occur on a daily basis. It is generally accepted that drones will play a key role in smart city environments in the near future. Drones will provide different types of public service, such as emergency response, disaster relief, and traffic control [2, 8]. In a smart city drones can be used to identify the location of accident, take pictures and notifying the emergency responders [9].

Least squares methods have been frequently used in regression analysis, where it leads approximate solutions and provides the best relationship between dependent and independent variables [10]. The method was originally discovered in 1795 by Carl Friedrich Gauss [11]. In addition to statistics, least squares methods stem from linear algebra and optimization theory [12]. They can be used for minimizing the sum of the squares in a set of equations to find the optimal solutions for various models in science, engineering, and several other fields [12]. This method can also benefits investors to understand the relationship between different variables. For example, in the world of crypto, an investor can understand the relationship between the coin share price and the earnings per coin [13]. In this paper, we are interested to use the concept of least squares solutions for traffic control and finding the best location of a drone observation system.

## 2. Mathematical Preliminaries

Orthogonal decomposition [14] is a mathematical method that can reduce the complexity of computationally intensive simulations such as computational fluid dynamics and structural analysis. For instance, it is typically used in flow and turbulence analysis to replace Navier-Stokes equations with models that are easier to solve. One of the most influential theorems in linear algebra is the orthogonal decomposition theorem, which is stated below [10]. The least squares method is based on the following theorem, which can be found in many linear algebra textbooks (see for example chapter 6 of [16]).

**Theorem 1. (Orthogonal decomposition theorem)** Let $W$ be a subspace of $R^n$. Then each **y** in $R^n$ can be written uniquely in the form

$$\mathbf{y} = \hat{\mathbf{y}} + \mathbf{z} \quad (1)$$

where $\hat{\mathbf{y}}$ is in $W$ and z is $W^\perp$. In fact, if $\{u_1, \ldots, u_p\}$ is any orthogonal basis of $W$, then

$$\hat{\mathbf{y}} = \frac{y.u_1}{u_1.u_1}\mathbf{u_1} + \ldots + \frac{y.u_p}{u_p.u_p}\mathbf{u_p} \quad (2)$$

and $\mathbf{z} = \mathbf{y} - \hat{\mathbf{y}}$.

The vector $\hat{y}$ in (1) is called the orthogonal projection of y into subspace $W$ and often is written as $\mathbf{proj}_W\mathbf{y}$. The proof of this theorem can be found on page 348 of the Linear Algebra book by David Lay [16]. This theorem lays the foundation for the method of least squares solution and the best approximation theorem, which is stated below.

**Theoren 2. (The Best Approximation Theorem)** Let **W** be a subspace of $R^n$. Let **y** be any vector in $R^n$, and let $\hat{y}$ be the orthogonal projection of y onto **W**, then $\hat{y}$ is the closest point in **W** to **y**, in the sense that

$$\|y - \hat{y}\| < \|y - v\| \quad (3)$$

for all **v** in **W** distinct from $\hat{y}$.

The vector $\hat{y}$ is called the best approximation to y by elements of W. The proof is provided in [16]. The best approximation theorem guarantees that the least squares solution of an inconsistent linear system is the minimum distance to all equations of the system. In our previous work [17] we presented the relationship between the symmedian point of a triangle and the least-squares solution of the corresponding linear system. Namely, the coordinates of the symmedian point are explicitly derived as the least squares solution of the linear system. This work was later extended to n-dimensional case, where the symmedian point and centroid of a petal tetrahedron can be determined using the least square method [18]. To explain this method, consider the following example of three lines forming a triangle on the xy-plane (see figure 1).

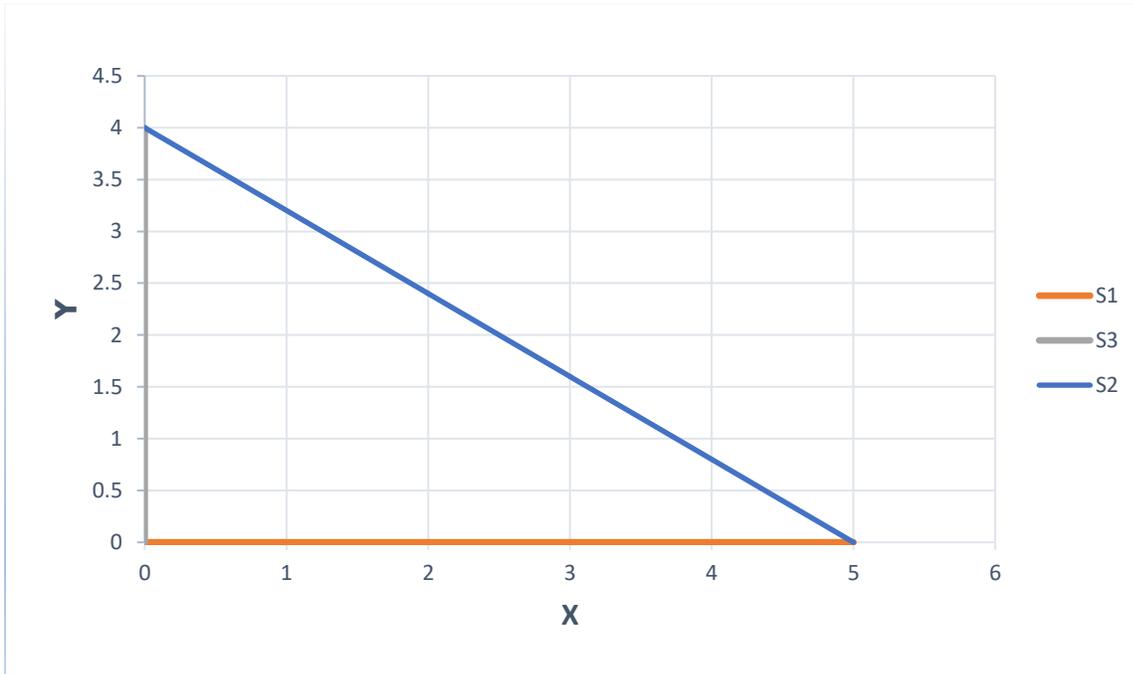

**Figure 1.** Graph of three lines presented in system (1)

**Example 1.** Consider the following system of three equation.

    S1: x = 0

    S2: 0.75x+y=4          (4)

    S3: y=0

As shown in Figure 1, system (4) is representing a triangle and there is no point that all three lines S1, S2 and S3 can pass through. Hence the system is inconsistent and does not have any solution. However, we can find the least squares solution, which has a minimized error according to equation (3) in Theorem 2. The main goal is to find the location $(\hat{x}, \hat{y})$ such that the sum of squared distance to each line S1, S2 and S3 is minimized. Specifically, the sum of the square system is given by $D(x,y)=d_1^2(x,y) + d_2^2(x,y) + d_3^2(x,y)$ and the objective is to minimize D by finding the min $D(x,y)=(\hat{x}, \hat{y})$. To do so, we need to take the following steps, which treats the problem in a general form of system

$$S1: a1x + b1y = c1$$
$$S2: a2x + b2y = c2 \quad (5)$$
$$S3: a3x + b3y = c3$$

**Step 1** Normalize the lines by dividing each equation by $\sqrt{ai^2 + bi^2}$ for $i = 1, 2, 3$.

$$A1 = \frac{a1}{\sqrt{a1^2+b1^2}}$$

$$B1 = \frac{b1}{\sqrt{a1^2+b1^2}}$$

$$C1 = \frac{c1}{\sqrt{a1^2+b1^2}}$$

$$\overline{S1}: 0x + 1y = 0$$
$$\overline{S2}: \frac{4}{\sqrt{41}} * (x) + (\frac{5}{\sqrt{41}})y = \frac{20}{\sqrt{41}}$$
$$\overline{S3}: 1x + 0y = 0$$

**Step 2** Find the least square solution using the formula:

$$(\hat{x}, \hat{y}) = (A^T.A)^{-1}(A^T.b)$$

For our example we have

$$A = \begin{pmatrix} 0 & 1 \\ \frac{4}{\sqrt{41}} & \frac{5}{\sqrt{41}} \\ 1 & 0 \end{pmatrix}, b = \begin{pmatrix} 0 \\ \frac{20}{\sqrt{41}} \\ 0 \end{pmatrix}, (A^T) = \begin{pmatrix} 0 & \frac{4}{\sqrt{41}} & 1 \\ 1 & \frac{5}{\sqrt{41}} & 0 \end{pmatrix}$$

$$(A^T.A) = \begin{pmatrix} 0 & \frac{4}{\sqrt{41}} & 1 \\ 1 & \frac{5}{\sqrt{41}} & 0 \end{pmatrix} * \begin{pmatrix} 0 & 1 \\ \frac{4}{\sqrt{41}} & \frac{5}{\sqrt{41}} \\ 1 & 0 \end{pmatrix} = \begin{pmatrix} \frac{57}{41} & \frac{20}{41} \\ \frac{20}{41} & \frac{66}{41} \end{pmatrix} = \frac{1}{41} * \begin{pmatrix} 57 & 20 \\ 20 & 66 \end{pmatrix}$$

$$(A^T.A)^{-1} = (\frac{1}{41}.\begin{pmatrix} 57 & 20 \\ 20 & 66 \end{pmatrix})^{-1} = 41 * \begin{pmatrix} 0.0196 & -0.0059 \\ -0.00599 & 0.0170 \end{pmatrix}$$

$$=\begin{pmatrix} 0.8049 & -0.2439 \\ -0.2439 & 0.6951 \end{pmatrix}$$

$$(A^T.b)=\begin{pmatrix} 0 & 0.75 & 1 \\ 1 & 1 & 0 \end{pmatrix}*\begin{pmatrix} 0 \\ 4 \\ 0 \end{pmatrix}=\begin{pmatrix} 1.9506 \\ 2.4384 \end{pmatrix} \text{ and } (A^T.A)^{-1}=\begin{pmatrix} 0.8049 & -0.2439 \\ -0.2439 & 0.6951 \end{pmatrix}$$

$$\begin{pmatrix} \hat{x} \\ \hat{y} \end{pmatrix} = (A^T.A)^{-1} *(A^T.b)=\begin{pmatrix} 0.8049 & -0.2439 \\ -0.2439 & 0.6951 \end{pmatrix}*\begin{pmatrix} 1.9506 \\ 2.4384 \end{pmatrix}=\begin{pmatrix} 0.9753 \\ 1.2192 \end{pmatrix}$$

So $(\hat{x}, \hat{y})$= (0.9753, 1.219) see figure 2 for the location of $(\hat{x}, \hat{y})$.

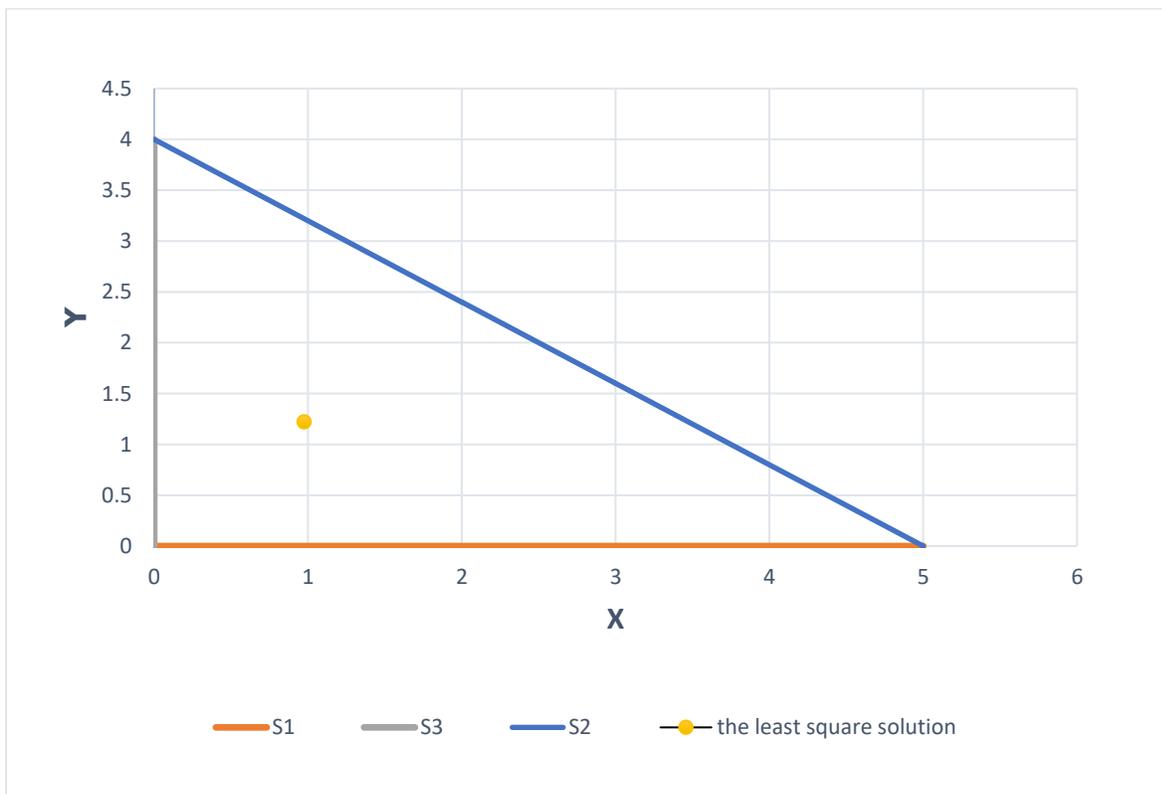

**Figure 2.** The least square solution of the linear system minimizes the square distances to the lines, which has been denoted in the figure.

## 3. Materials and Method

In this paper we are extending this concept to the case multiple lines creating a polygon. Then we apply this method to the real-world example of traffic control in smart cities. Note that each equation must be normalized. An alternative method is to use the weighted least squares method which normalizes the equations by multiplying a diagonal matrix W to both sides of equation AX = b. Kansas City's traffic data can be obtained from different online resources such as TomTom [19].

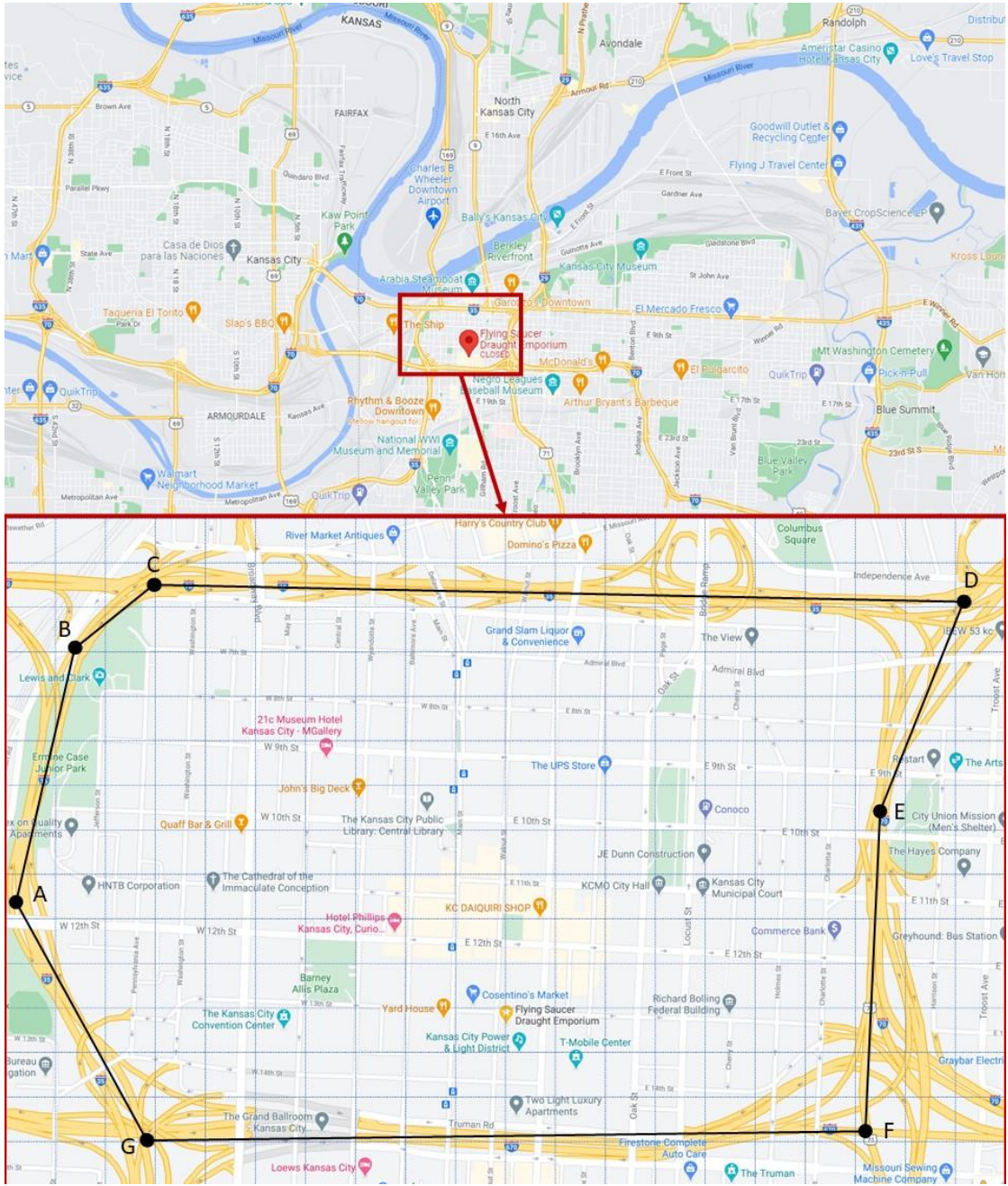

**Figure 3.** Weighted least squares method applied to Power and Light District of Kansas City, Missouri. Top. Map of Kansas City, Missouri, Bottom. Power and Light District.

Here we would like to identify the best location of a drone station in Kansas City Power and Light District based on the percent of accidents and traffic congestions for each of the highways surrounding the region (See Figure 3).

We take the following steps to find the optimal location of the drone station.

a) Determine the approximate coordinates of points A-G (see figure 3)
b) Find the equations of lines AB – GA

**Important note:** After obtaining the equations of the lines we need to normalize them. For instance, if $ax + by = c$ is one of the lines you must change it to $(a/\sqrt{a^2 + b^2})x + (b/\sqrt{a^2 + b^2})y = (c/\sqrt{a^2 + b^2})$

c) Make the suitable for system of equations AX=b corresponding to lines AB – GA
d) Determine the diagonal weight matrix W using the morning rush hour (RH) percent accident/congestion data shown in Table 1.

**Table 1)** Summary data of percent congestions for each highway segment during morning and afternoon Rush Hours (RH) in Kansas City Power and Light District.

| Route | AB | BC | CD | DE | EF | FG | GA |
|---|---|---|---|---|---|---|---|
| **Morning RH** | 2% | 9% | 3% | 25% | 45% | 15% | 1% |
| **Afternoon RH** | 17% | 10% | 30% | 15% | 6% | 10% | 12% |

e) Using parts (c) and (d) determine the best morning location of the drone station (i.e., the location that has the minimum sum of weighted squared distances). Specifically, find the location of $(\hat{x}, \hat{y})$ using weighed lease squares by doing the following computations:

- multiply both sides of AX=b by W to get WAX=Wb
- multiply both sides of WAX=Wb by $(WA)^T$ to get $(WA)^T$ WAX=$(WA)^T$ Wb
- multiply both sides of $(WA)^T$ WAX=$(WA)^T$ Wb by $((WA)^T WA)^{-1}$ to get the least square solution $(\hat{x}, \hat{y})$.

f) Determine the best morning location of the drone station using rush hour percent accident/congestion data.

## 4. Results

Using the methodology outlined in the previous section, we find the optimal location of drone station according to the morning and afternoon rush hour data. As detailed below, the computation of weighted least squares we will result in two stations, one for the morning and therefore the afternoon. As shown in Figure 4, the equations of lines and their intersections can be determined using basic arithmetic. Next, we apply the following steps to the corresponding lines to find the least squares solution without rush hour data

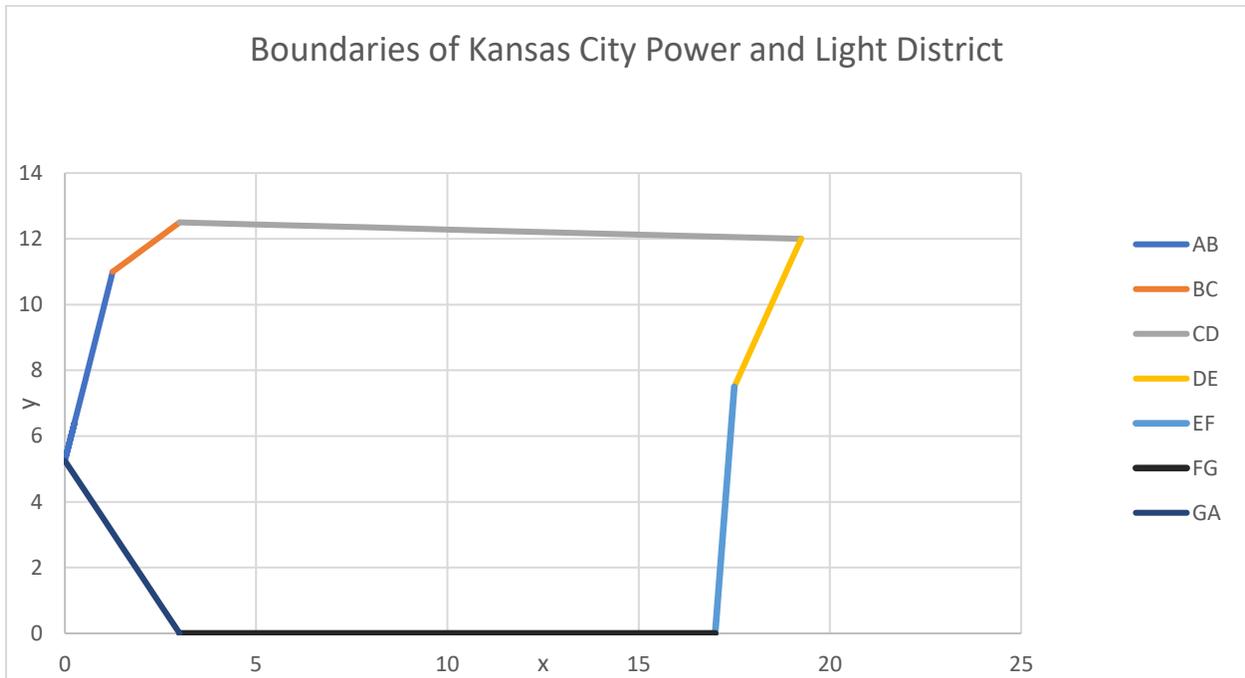

**Figure 4.** Each line represents a specific segment of a highway surrounding the Kansas City Power and Light District

**Table 2)** Summary data of intersections and equations of lines

| Interaction points | Equations of Lines | Parametrized Equations |
|---|---|---|
| A (0,5.25) | AB=S1: a1x+b1y = C1 | S1:  -4.6x+ y =5.25 |
| B (1.25,11) | BC= S2: a2x+b2y = C2 | S2   y-0.857x = 9.928 |
| C (3,12.5) | CD= S3: a3x+b3y = C3 | S3   y + 0.03076 x = 12.5923 |
| D (19.25,12) | DE= S4: a4x+b4y = C4 | S4: 2.5714 x – y = 37.5 |
| E (17.5,7.5) | EF= S5: a5x+b5y = C5 | S5: 15x -y = 255 |
| F (17,0) | FG= S6: a6x+b6y = C6 | S6:  y =0 |
| G (3,0) | GA= S2: a7x+b7y = C7 | S7: 1.75x+y =5.25 |

**Step 1)** Normalize the line by dividing each equation by $\sqrt{ai^2 + bi^2}$

$$A1 = \frac{a1}{\sqrt{a1^2+b1^2}},\ B1 = \frac{b1}{\sqrt{a1^2+b1^2}},\ C1 = \frac{c1}{\sqrt{a1^2+b1^2}}$$

$$S1: \frac{-4.6}{\sqrt{(-4.6)^2 + 1}} * (x) + (\frac{1}{\sqrt{(-4.6)^2 + 1}})y = \frac{5.25}{\sqrt{(-4.6)^2 + 1}}$$

$$S1: \ -0.977\,x + 0.212\,y = 1.115$$

$$S2: \frac{-0.857}{\sqrt{(-0.857)^2 + 1}} * (x) + \left(\frac{1}{\sqrt{(-0.857)^2 + 1}}\right) y = \frac{9.928}{\sqrt{(-0.857)^2 + 1}}$$

$$S2: \ -0.652\,x + 0.761\,y = 7.56$$

$$S3: \frac{0.03076}{\sqrt{(0.03076)^2 + 1}} * (x) + \left(\frac{1}{\sqrt{(0.03076)^2 + 1}}\right) y = \frac{12.5923}{\sqrt{(0.03076)^2 + 1}}$$

$$S3: \ 0.0307\,x + 1\,y = 12.586$$

$$S4: \frac{2.5714}{\sqrt{(2.5714)^2 + 1}} * (x) - \left(\frac{1}{\sqrt{(2.5714)^2 + 1}}\right) y = \frac{37.5}{\sqrt{(2.5714)^2 + 1}}$$

$$S4 \quad 0.932\,x - 0.362\,y = 13.59$$

$$S5: \frac{15}{\sqrt{(15)^2 + 1}} * (x) - \left(\frac{1}{\sqrt{(15)^2 + 1}}\right) y = \frac{255}{\sqrt{(15)^2 + 1}}$$

$$S5: \ 0.977x - 0.0665\,y = 1.115$$

$$S6: \frac{0}{\sqrt{(0)^2 + 1}} * (x) + \left(\frac{1}{\sqrt{(0)^2 + 1}}\right) y = \frac{0}{\sqrt{(0)^2 + 1}}$$

$$S6: 0\,x + y = 0$$

$$S7: \frac{1.75}{\sqrt{(1.75)^2 + 1}} * (x) + \left(\frac{1}{\sqrt{(1.75)^2 + 1}}\right) y = \frac{5.25}{\sqrt{(1.75)^2 + 1}}$$

$$S7: \ 0.868\,x + 0.4961\,y = 2.604$$

**Outputs of step 1)**

S1: - 0.977 x +0.212 y = 1.115
S2: - 0.652 x +0.761 y = 7.56
S3: 0.0307 x + 1 y = 12.586
S4  0.932 x - 0.362 y = 13.59
S5: 0.977x -0.0665 y = 1.115
S6: 0x   + y =0
S7: 0.868 x +0.4961 y = 2.604

The set of equations S1-S7 is equivalent to the linear system $AX = b$, where

$$A = \begin{pmatrix} -0.977 & 0.212 \\ -0.652 & 0.761 \\ 0.0307 & 1 \\ 0.932 & -0.362 \\ 0.977 & -0.0665 \\ 0 & 1 \\ 0.868 & 0.4961 \end{pmatrix} \quad \text{and} \quad b = \begin{pmatrix} 1.115 \\ 7.56 \\ 12.586 \\ 13.59 \\ 1.115 \\ 0 \\ 2.604 \end{pmatrix}$$

**Step 2** Find the least square solution using the formula:

$$(\hat{x}, \hat{y}) = (A^T.A)^{-1} (A^T.b)$$

$$(A^T) = \begin{pmatrix} -\frac{977}{1000} & -\frac{163}{250} & \frac{307}{10000} & \frac{233}{250} & \frac{977}{1000} & 0 & \frac{217}{250} \\ \frac{53}{250} & \frac{76181}{1000} & 1 & -\frac{181}{500} & -\frac{133}{2000} & 1 & \frac{4961}{10000} \end{pmatrix}$$

$$(A^T.A) = \begin{pmatrix} -\frac{977}{1000} & -\frac{29}{50} & \frac{3}{200} & \frac{461}{500} & \frac{539}{1000} & 0 & \frac{171}{200} \\ \frac{53}{250} & \frac{81}{100} & 1 & -\frac{2}{5} & -\frac{179}{1000} & 1 & \frac{1}{2} \end{pmatrix} * \begin{pmatrix} -0.977 & 0.212 \\ -0.58 & 0.81 \\ 0.015 & 1 \\ 0.922 & -0.4 \\ 5.39 & -0.179 \\ 0 & 1 \\ 0.87 & 0.5 \end{pmatrix} =$$

$$\begin{pmatrix} \frac{395715249}{100000000} & -\frac{6443357}{10000000} \\ -\frac{6443357}{10000000} & \frac{150282323}{50000000} \end{pmatrix}$$

$$(A^T.A)^{-1} = \begin{pmatrix} \frac{15028232300000000}{57393164394770977} & \frac{3221678500000000}{57393164394770977} \\ \frac{3221678500000000}{57393164394770977} & \frac{19785762450000000}{57393164394770977} \end{pmatrix}$$

$$(A^T.b) = \begin{pmatrix} -\frac{977}{1000} & -\frac{29}{50} & \frac{3}{200} & \frac{461}{500} & \frac{539}{1000} & 0 & \frac{171}{200} \\ \frac{53}{250} & \frac{81}{100} & 1 & -\frac{2}{5} & -\frac{179}{1000} & 1 & \frac{1}{2} \end{pmatrix} * \begin{pmatrix} 1.115 \\ 7.56 \\ 12.586 \\ 13.59 \\ 1.115 \\ 0 \\ 2.604 \end{pmatrix} = \begin{pmatrix} \frac{51917111}{5000000} \\ \frac{148736569}{10000000} \end{pmatrix}$$

$$(A^T.A)^{-1} * (A^T.b) = \begin{pmatrix} \frac{51917111}{5000000} \\ \frac{148736569}{10000000} \end{pmatrix} * \begin{pmatrix} \frac{15028232300000000}{57393164394770977} & \frac{3221678500000000}{57393164394770977} \\ \frac{3221678500000000}{57393164394770977} & \frac{19785762450000000}{57393164394770977} \end{pmatrix}$$

$$= \begin{pmatrix} \dfrac{203962621541683710}{57393164394770977} \\ \dfrac{327738690244366105}{57393164394770977} \end{pmatrix} = \begin{pmatrix} 3.55 \\ 5.71 \end{pmatrix}$$

Hence, the least squares solution is given by $\hat{x} = 3.55$, $\hat{y} = 5.71$ as shown in Figure 5.

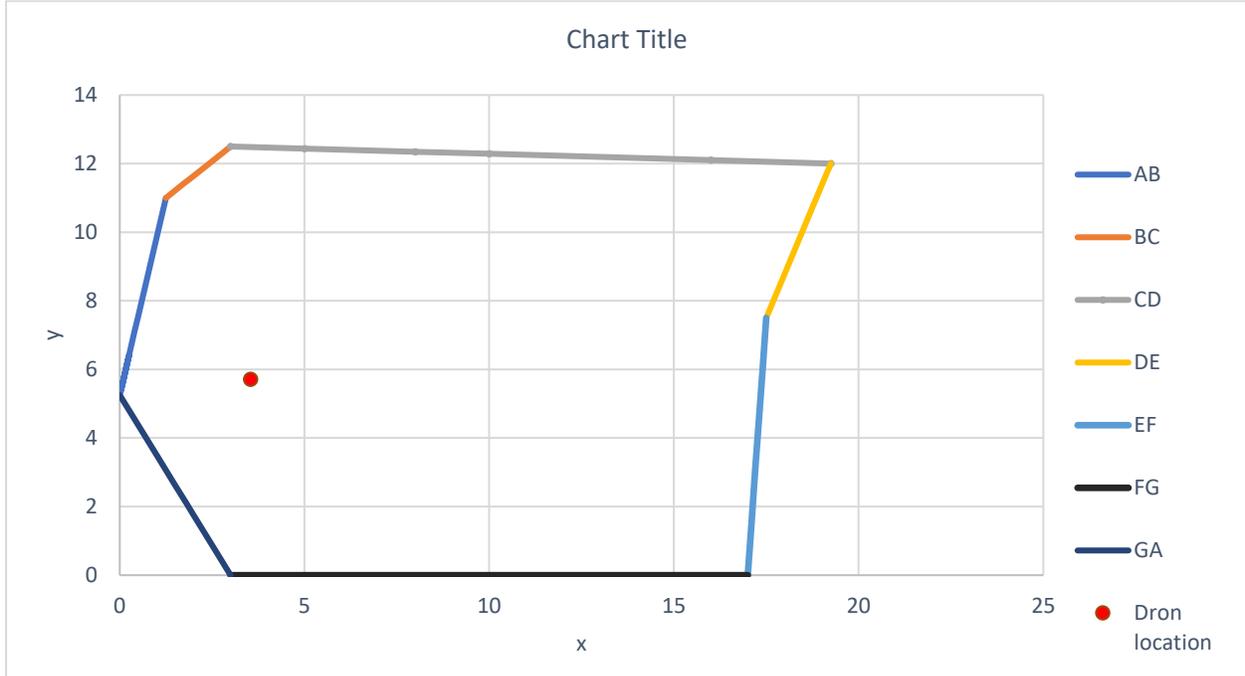

**Figure 5.** Least squares solution without rush hour data. We get that $\hat{x} = 3.55$ and $\hat{y} = 5.71$ as shown in the map.

Using the weighted least the squares approach, we can apply the morning and afternoon rush hour data to update the optimal location of drone station. There is an extra step that should be followed: The linear system $AX = b$ should be multiplied by weight matrix $W_M$ and $W_A$ which contain the morning and afternoon rush hour data, respectively. We have

$$W_M = \begin{bmatrix} 0.2 & 0 & 0 & 0 & 0 & 0 & 0 \\ 0 & 0.9 & 0 & 0 & 0 & 0 & 0 \\ 0 & 0 & 0.3 & 0 & 0 & 0 & 0 \\ 0 & 0 & 0 & 0.25 & 0 & 0 & 0 \\ 0 & 0 & 0 & 0 & 0.45 & 0 & 0 \\ 0 & 0 & 0 & 0 & 0 & 0.15 & 0 \\ 0 & 0 & 0 & 0 & 0 & 0 & 0.1 \end{bmatrix}$$

$$A=\begin{pmatrix} -0.977 & 0.212 \\ -0.652 & 0.761 \\ 0.0307 & 1 \\ 0.932 & -0.362 \\ 0.977 & -0.0665 \\ 0 & 1 \\ 0.868 & 0.4961 \end{pmatrix} \quad b=\begin{pmatrix} 1.115 \\ 7.56 \\ 12.58 \\ 13.59 \\ 1.115 \\ 0 \\ 2.604 \end{pmatrix}$$

$$(W_M*A)=\begin{bmatrix} 0.2 & 0 & 0 & 0 & 0 & 0 & 0 \\ 0 & 0.9 & 0 & 0 & 0 & 0 & 0 \\ 0 & 0 & 0.3 & 0 & 0 & 0 & 0 \\ 0 & 0 & 0 & 0.25 & 0 & 0 & 0 \\ 0 & 0 & 0 & 0 & 0.45 & 0 & 0 \\ 0 & 0 & 0 & 0 & 0 & 0.15 & 0 \\ 0 & 0 & 0 & 0 & 0 & 0 & 0.1 \end{bmatrix} \begin{pmatrix} -0.977 & 0.212 \\ -0.652 & 0.761 \\ 0.0307 & 1 \\ 0.932 & -0.362 \\ 0.977 & -0.0665 \\ 0 & 1 \\ 0.868 & 0.4961 \end{pmatrix} =$$

$$\begin{pmatrix} -0.1954 & 0.0424 \\ -0.5220 & 0.7290 \\ 0.0060 & 0.3 \\ 0.2325 & -0.0900 \\ 0.4455 & -0.0297 \\ 0 & 0.15 \\ 0.0860 & 0.0490 \end{pmatrix}$$

$$(W_M*b)=\begin{bmatrix} 0.2 & 0 & 0 & 0 & 0 & 0 & 0 \\ 0 & 0.9 & 0 & 0 & 0 & 0 & 0 \\ 0 & 0 & 0.3 & 0 & 0 & 0 & 0 \\ 0 & 0 & 0 & 0.25 & 0 & 0 & 0 \\ 0 & 0 & 0 & 0 & 0.45 & 0 & 0 \\ 0 & 0 & 0 & 0 & 0 & 0.15 & 0 \\ 0 & 0 & 0 & 0 & 0 & 0 & 0.1 \end{bmatrix} \begin{pmatrix} 1.115 \\ 7.56 \\ 12.3 \\ 13.6 \\ 16.9 \\ 0 \\ 2.61 \end{pmatrix} = \begin{pmatrix} 0.223 \\ 0.819 \\ 3.69 \\ 3.4 \\ 7.6050 \\ 0 \\ 0.261 \end{pmatrix}$$

multiply both sides of WAX=Wb by $(WA)^T$ to get $(WA)^T WAX = (WA)^T Wb$

$(Wm*A)^T =$
$$\begin{pmatrix} -0.1954 & -0.522 & 0.006 & 0.2325 & 0.4455 & 0 & 0.086 \\ 0.0424 & 0.729 & 0.3 & 0.09 & -0.0297 & 0.15 & 0.049 \end{pmatrix}$$

$$(Wm*A)^T *(Wm*A) = \begin{pmatrix} -0.1954 & -0.522 & 0.006 & 0.2325 & 0.4455 & 0 & 0.086 \\ 0.0424 & 0.729 & 0.3 & 0.09 & -0.0297 & 0.15 & 0.049 \end{pmatrix} *$$

$$\begin{pmatrix} -0.1954 & 0.0424 \\ -0.5220 & 0.7290 \\ 0.0060 & 0.3 \\ 0.2325 & -0.0900 \\ 0.4455 & -0.0297 \\ 0 & 0.15 \\ 0.0860 & 0.0490 \end{pmatrix}$$

$$((Wm*A).(Wm*A)^T)^{-1} = \begin{pmatrix} 3.2675 & 2.0733 \\ 2.0733 & 2.8374 \end{pmatrix}$$

$$(Wm*A)^T.(W_m*b) =$$

$$\begin{pmatrix} -0.1954 & -0.522 & 0.006 & 0.2325 & 0.4455 & 0 & 0.086 \\ 0.0424 & 0.729 & 0.3 & 0.09 & -0.0297 & 0.15 & 0.049 \end{pmatrix} * \begin{pmatrix} 0.223 \\ 0.819 \\ 3.69 \\ 3.4 \\ 7.6050 \\ 0 \\ 0.261 \end{pmatrix} =$$

$$\begin{pmatrix} 3.752 \\ 1.1944 \end{pmatrix}$$

Now, multiply both sides of $(WA)^T WAX = (WA)^T Wb$ by $((WA)^T WA)^{-1}$ to get the weighted least squares solution

$$X = \begin{bmatrix} x \\ y \end{bmatrix} = ((WA)^T WA)^{-1} Wb \ ((WA)^T * WA)^{-1} =$$

$$(Wm*A)^T.(W_m*b)*((Wm*A).(Wm*A)^T)^{-1} =$$

$$\begin{pmatrix} 3.752 \\ 1.1944 \end{pmatrix} * = \begin{pmatrix} 3.2675 & 2.0733 \\ 2.0733 & 2.8374 \end{pmatrix}$$

$$= \begin{pmatrix} 14.7361 \\ 11.1682 \end{pmatrix}$$

Hence, the weighted least squares solution is given by $\hat{x} = 3.3$, $\hat{y} = 3.7$ as shown in Figure 6.

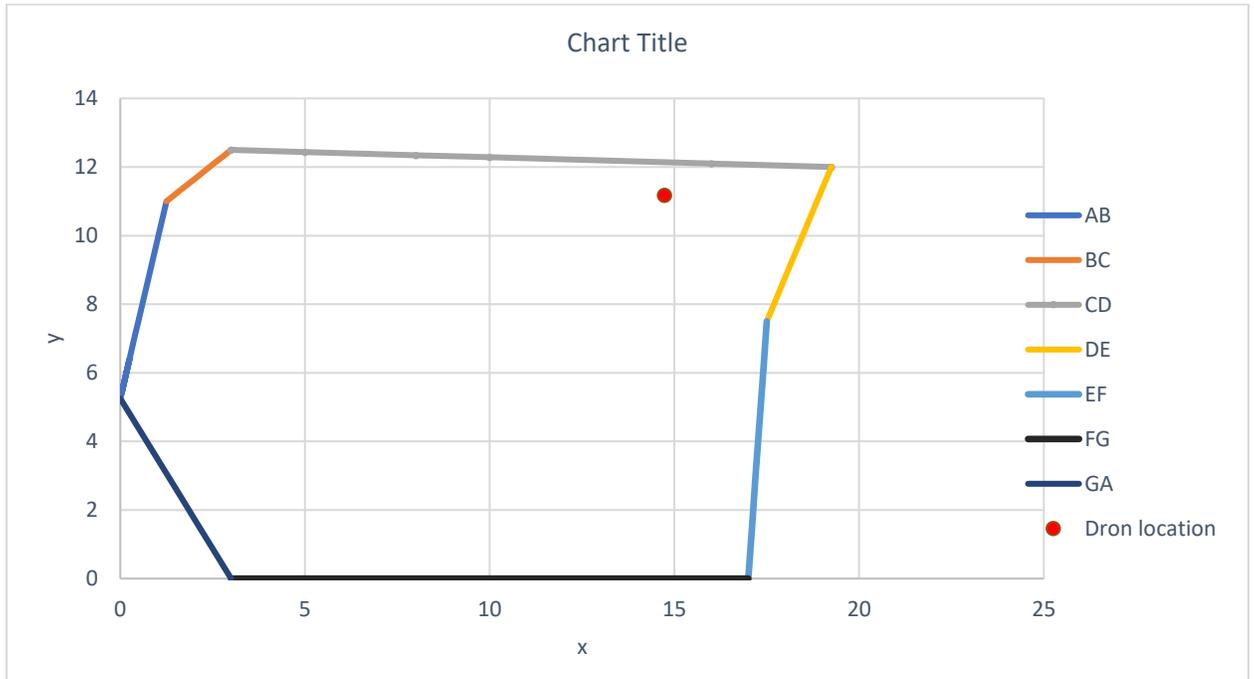

**Figure 6.** Least squares solution morning rush hour data. We get that $\hat{x} = 3.3$, $\hat{y} = 3.7$ as shown in the map.

Similarly, we can apply the weight matrix $W_A$ for the afternoon

$$W_A = \begin{bmatrix} 0.17 & 0 & 0 & 0 & 0 & 0 & 0 \\ 0 & 0.10 & 0 & 0 & 0 & 0 & 0 \\ 0 & 0 & 0.30 & 0 & 0 & 0 & 0 \\ 0 & 0 & 0 & 0.15 & 0 & 0 & 0 \\ 0 & 0 & 0 & 0 & 0.6 & 0 & 0 \\ 0 & 0 & 0 & 0 & 0 & 0.10 & 0 \\ 0 & 0 & 0 & 0 & 0 & 0 & 0.12 \end{bmatrix}$$

Following the same procedure, we get the optimized drone location for the afternoon should be

$$A=\begin{pmatrix} -0.977 & 0.212 \\ -0.652 & 0.761 \\ 0.0307 & 1 \\ 0.932 & -0.362 \\ 0.977 & -0.0665 \\ 0 & 1 \\ 0.868 & 0.4961 \end{pmatrix} \quad b=\begin{pmatrix} 1.115 \\ 7.56 \\ 12.58 \\ 13.59 \\ 1.115 \\ 0 \\ 2.604 \end{pmatrix}$$

$$(W_M * A) = \begin{pmatrix} -0.1661 & 0.036 \\ -0.0580 & 0.081 \\ 0.0060 & 0.3 \\ 0.1395 & -0.054 \\ 0.5940 & -0.0396 \\ 0 & 0.1 \\ 0.1032 & 0.0588 \end{pmatrix}$$

$$(Wa * A)^T = \begin{pmatrix} -0.1661 & -0.058 & 0.006 & 0.1395 & 0.594 & 0 & 0.1032 \\ 0.0360 & 0.081 & 0.3 & 0.05 & -0.0396 & 1 & 0.0588 \end{pmatrix}$$

$$(Wa * A)^T * (Wa * A) = \begin{pmatrix} 0.4175 & -0.0177 \\ -0.0177 & 0.1158 \end{pmatrix}$$

$$((Wa * A)^T (Wa * A))^{-1} = \begin{pmatrix} 2.4108 & 0.3679 \\ 0.3679 & 8.6916 \end{pmatrix}$$

$$(Wa * b) = \begin{pmatrix} 0.1896 \\ 0.091 \\ 3.69 \\ 2.04 \\ 10.14 \\ 0 \\ 0.3132 \end{pmatrix}$$

$$(Wa * A)^T \cdot (Wa * b) = \begin{pmatrix} 6.5247 \\ 0.6279 \end{pmatrix}$$

$$((Wa * A)^T (Wa * A))^{-1} * (Wa * A)^T \cdot (Wa * b) = \begin{pmatrix} 2.4108 & 0.3679 \\ 0.3679 & 8.6916 \end{pmatrix} * \begin{pmatrix} 6.5247 \\ 0.6279 \end{pmatrix} = \begin{pmatrix} 15.96 \\ 7.8579 \end{pmatrix}$$

$$\begin{bmatrix} x \\ y \end{bmatrix} = \begin{bmatrix} 15.9608 \\ 7.857 \end{bmatrix}.$$

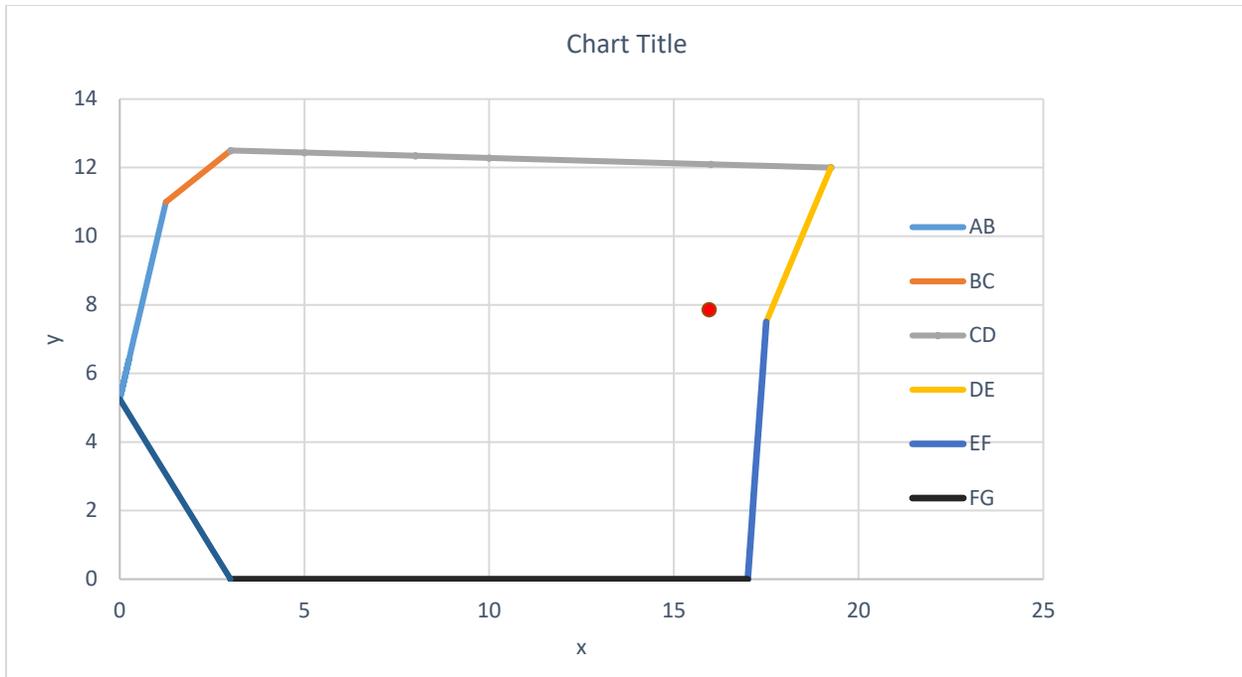

The main purpose of this paper was to highlight the applicability of least squares method in smart cities. Using the example of Kansas City, Missouri, we showed that the optimal location of a portable drone station can be updated according to the morning and afternoon rush hour data. It should be noted that the computations presented in the paper are mainly for the purpose of understanding how the weighted least squares method works. Otherwise, such computations can be carried out and updated using GPS and AI systems designed for smart cities [20].

[20] Ullah, Zaib, Fadi Al-Turjman, Leonardo Mostarda, and Roberto Gagliardi. "Applications of artificial intelligence and machine learning in smart cities." Computer Communications 154 (2020): 313-323.18